\documentclass[a4paper,12pt]{article}

\usepackage[latin1]{inputenc}
\usepackage{amsmath,amscd, dsfont}
\usepackage{amssymb,epsf,graphics,amsmath}
\usepackage{graphicx}
\usepackage{ifthen}
\usepackage{makeidx}
\usepackage{mathrsfs}
\usepackage{pstcol}
\usepackage[T1]{fontenc}

\newtheorem{theo}{Theorem}[section]
\newtheorem{prop}[theo]{Proposition}

\newtheorem{lem}[theo]{Lemma}

\newenvironment{proof}{{\bf Proof }}{\hfill $\Box$}

\newcommand{\C}{\mathbb{C}}
\newcommand{\N}{\mathbb{N}}
\newcommand{\R}{\mathbb{R}}

\newcommand{\Bb}{\mathcal{B}}

\newcommand{\Hh}{\mathcal{H}}

\newcommand{\Pp}{\mathcal{P}}

\newcommand{\Tr}{\mathrm{Tr}}

\begin{document}


\title{\sc Low density limit and the quantum Langevin equation for the heat bath}
\author{Ameur DHAHRI \\\vspace{-2mm}\scriptsize Ceremade, UMR CNRS 7534\\\vspace{-2mm}\scriptsize Universit\'e Paris Dauphine\\\vspace{-2mm}\scriptsize place De Lattre de Tassigny\\\vspace{-2mm}\scriptsize 75775 Paris C\'edex 16\\\vspace{-2mm}\scriptsize France \\\vspace{-2mm}\scriptsize e-mail: dhahri@ceremade.dauphine.fr}

\date{}

\maketitle
\begin{abstract}
  We consider a repeated quantum interaction model describing a small
  system $\Hh_S$ in interaction with each one of the identical copies of the chain $\bigotimes_{\N^*}\C^{n+1}$, modeling a heat bath, one after another
  during the same short time intervals $[0,h]$. We suppose that the
  repeated quantum interaction Hamiltonian is split in two parts: a
  free part and an interaction part with time scale of order $h$.
  After giving the GNS representation, we establish the connection
  between the time scale $h$ and the classical low density limit. We
  introduce a chemical potential $\mu$ related to the time $h$ as
  follows: $h^2=e^{\beta\mu}$. We further prove that the solution of
  the associated discrete evolution equation converges strongly, when
  $h$ tends to 0, to the unitary solution of a quantum Langevin
  equation directed by Poisson processes.
\end{abstract}
\section{Introduction}
In the quantum theory of open systems, two different approaches have
usually been considered by physicists as well as mathematicians: The
Hamiltonian and Markovian approaches.

The first approach consists in giving a full Hamiltonian description
for the interactions of a quantum system with a quantum field
(reservoir, heat bath...) and studying the ergodic properties of the
associated dynamical system.

The second approach consists of giving up the idea of modeling the
quantum field and concentrating on the effective dynamics of the
quantum system. The dynamics are then described by a Lindblad generator,
which dilates a quantum Langevin equation (or quantum stochastic differential equation (cf [P])).

It is worthwhile to note that the quantum Langevin equation (or Lindblad
generator) associated to the combined system can be derived from its
Hamiltonian description by the classical weak coupling and low density
limits (cf [AFL], [APV], [Dav], [D1], [Pe]...).

Recently, Attal and Pautrat describe the interaction between a quantum
system and a quantum field by a {\it repeated quantum interaction
  model} (cf [AtP]): The exterior system is modeled by an
infinite chain of identical copies $\mathcal{H}$ ($\mathcal{H}$ is an
Hilbert space) and the interaction is described as follows: The small
system interacts with each one of the identical pieces of the exterior
system one after another during the same time intervals $[0,h]$. They
prove, in the continuous limit ($h$ tends to 0), that this discrete
description of the combined system gives rise to a quantum Langevin
equation.

In [AtP], we observe three time scales, which appear in a repeated
quantum interaction Hamiltonian $H=H(h)$, with respectively order 1, $\sqrt{h}$ and $h$. 

In [AtJ], the authors have studied the time scale of order $\sqrt{h}$. They prove that, in the continuous limit, we get  a quantum diffusion equation where new noises, called {\it thermal noises}, appear. This normalization is used for modeling some physical systems (cf [D1], [D2]).

In this paper we prove that the time scale of order $h$ corresponds to the low density limit. For this purpose, we consider a repeated interaction model associated to a small system  $\Hh_S$, with Hamiltonian $H_S$, which interacts with a chain $\bigotimes_{\N^*}\C^{n+1}$ of a heat bath, so that the Hamiltonian of each piece of the chain is the operator $H_R$ on $\C^{n+1}$. The associated repeated quantum interaction Hamiltonian is defined on $\Hh_S\otimes\C^{n+1}$ by
\begin{eqnarray}\label{Hamil}
H=H_S\otimes I+I\otimes H_R+\frac{1}{h}\sum_{i,j=1}^nD_{ij}\otimes a^i_j,
\end{eqnarray}
where the $D_{ij}$ are the interaction operators associated to the small system and the $a^i_j$ are the discrete quantum noises in $\Bb(\C^{n+1})$, the algebra of all bounded operators on $\C^{n+1}$. The thermodynamical equilibrium state of a one copy of the atom chain is defined by the density matrix
$$\rho_\beta=\frac{e^{-\beta(H_R-\mu N)}}{\Tr(e^{-\beta(H_R-\mu N)})}\;,$$
where $N$ is the discrete number operator defined on $\C^{n+1}$ and $\mu$ is a chemical potential ($\mu<0$).

The length of the time interaction $h$ between the small system and each piece of the heat bath is supposed to be related to the chemical potential by the relation
$$h^2=e^{\beta\mu}.$$ 
Obviously, $h$ tends to $0$ if and only if $\mu$ tends to $-\infty$.

After giving the GNS representation and taking into account the above assumption, we prove that in the continuous time limit we get a quantum stochastic differential equation directed by Poisson processes.

This paper is organized as follows. In Section 2 we introduce the discrete model which presents the repeated interaction model describing a small system in interaction with a heat bath. Also, we give a description of the GNS representation of the pair $(\C^{n+1},\,\rho_\beta)$. In Section 3 we describe the tools used to obtain the continuous limit: Guichardet interpretation of a Fock space, quantum noises and quantum Langevin equations. Finally, in Section 4 we prove that the discrete solution of the associated discrete evolution equation with repeated quantum interaction Hamiltonian, given by (\ref{Hamil}), converges strongly to the unitary solution of a quantum Langevin equation. In this equation only Poisson processes appear in its noise part.
\section{The discrete model}
In this section we start by describing the discrete atom chain modeling an exterior system (reservoir, heat bath...). We further give the repeated quantum interaction model, which is the object of our study. In the last part we describe the associated GNS representation.          
\subsection{The atom chain}\label{subsection}
Let us give a brief description of the algebraic structure of the atom
chain. We refer the interested reader to [At2] for more details. Let
$\Hh$ be a Hilbert space, where we fix an orthonormal basis
$\{e_k,\;k\in J=I\cup\{0\}\}$. The vector $e_0=\Omega$ defines the
vacuum state.  Now, consider the {\it atom chain} $T\Phi=\otimes_{\N^*}\Hh$ defined with respect to the stabilizing sequence
$(\Omega)_{n\in\mathbb{N}^*}$ and denote by $\Pp_{\N^*,J}$ the set of finite subsets $\{(n_1,i_1),...,(n_k,i_k)\}$ of $\N^*\times J$ such that $n_i\neq n_j$ for all $i\neq j$. Then, an orthonormal basis of $T\Phi$ is given by the family
$$\{e_\sigma,\;\sigma\in\Pp_{\N^*,J}\},$$
where $e_\sigma,\,\sigma=\{(n_1,i_1),...,(n_k,i_k)\}$, is the infinite tensor product of elements of the basis $ \{e_k,\;k\in J=I\cup\{0\}\}$ such that $e_{i_{m}}\;(1\leq m\leq k)$ appears in the $n_{i_{m}}$-th copy of $\Hh$ and $\Omega$ appears in the other copies of $\Hh$ in the tensor product $\otimes_{\N^*}\Hh$.

Let $\{a^i_j,\; i,j\in J\}$ be the basis of $\Bb(\Hh)$ defined by 
$$ a^i_j e_k=\delta_{ik}e_j.$$
Denote by $a^i_j(k)$ the operators on $T\Phi$, which act as $a^i_j$ on the $k$-th copy of $\Hh$ in the atom chain $\otimes_{\N^*}\Hh$ and the identity elsewhere. The operators $a^i_j(k)$ are called {\it discrete quantum noises} and they act on elements of the basis $\{e_\sigma,\;\sigma\in\Pp_{\N^*,J}\}$ as follows
\begin{eqnarray*}
a^i_j(k)e_\sigma&=&\mathds{1}_{(k,i)\in\sigma}e_{\sigma\setminus(k,i)\cup(k,j)},\; \mbox{ for all }\,i\neq0 \mbox{ and }\, j\neq0,\\
a^i_0(k)e_\sigma&=&\mathds{1}_{(k,i)\in\sigma}e_{\sigma\setminus(k,i)},\,\mbox{ for all }\,i\neq0,\\
a^0_j(k)e_\sigma&=&\mathds{1}_{\{(k,i)\notin \sigma,\forall i\in J\}}e_{\sigma\cup(k,j)},\,\mbox{ for all }\,j\neq0,\\
a^0_0(k)e_\sigma&=&\mathds{1}_{\{(k,i)\notin\sigma,\forall i\in I\}}e_{\sigma}.
\end{eqnarray*}
\subsection{Small system in interaction with a heat bath}
Now, we consider a small system described by a Hilbert space $\Hh_S$
in interaction with a heat bath modeled by the atom chain
$\bigotimes_{\N^*}\C^{n+1}$, where $\Bb=\{e_0,e_1,...,e_n\}$ is an
orthonormal basis of $\Hh=\C^{n+1}$ and $e_0=\Omega$ is the vacuum
state. The interaction between the two systems is described as
follows: the small system interacts with each of the identical copies $\C^{n+1}$ of the heat bath one after another during the consecutive short time intervals $[nh,(n+1)h]$. Therefore, the total interaction between the small system and the chain of identical pieces is described by the Hilbert space $\Hh_S\otimes\bigotimes_{\N^*}\C^{n+1}$. 

Consider the orthonormal basis $\{a^i_j,\;0\leq i,j\leq n\}$ of $\Bb(\C^{n+1})$ where
$$a^i_je_k=\delta_{ik}e_j.$$
The respective Hamiltonians of the small system and one piece of the heat bath are given by a self-adjoint operator $H_S$ defined on $\Hh_S$ and the operator $H_R$ defined on $\C^{n+1}$ by
$$H_R=\sum_{i=0}^n\gamma_i a_i^0a_0^i,$$
where $\gamma_i$ are real numbers.

The full Hamiltonian of the small system interacting with one piece is the self-adjoint operator $H$ defined on $\Hh_S\otimes\C^{n+1}$ by
\begin{eqnarray}\label{scale}
H=H_S\otimes I+I\otimes H_R+\frac{1}{h}\sum_{i,j=1}^n\,D_{ij}\otimes a^i_j,
\end{eqnarray}
where $D_{ij}=(D_{ji})^*$. Note that the operators $a^i_j$ describe the transition of a non-empty state of one piece to another non-empty state, where the total number of particles is preserved.

The associated unitary evolution during the time interval $[0,h]$ is the operator
$$U=e^{-ihH}.$$

Denote by $\C^{n+1}_k$ the $k$-copy of $\C^{n+1}$ in the chain $\bigotimes_{\N^*}\C^{n+1}$. Then we define the operator $U_k$ on $\Hh_S\otimes\bigotimes_{\N^*}\C^{n+1}$ by
\begin{eqnarray*}
U_k=
\left\{
\begin{array}{lcc}
U&&\mbox {on }\,\Hh_S\otimes\C^{n+1}_k\\
I&& \mbox{ elsewhere }.
\end{array}
\right.
\end{eqnarray*}
Hence, the discrete evolution equation, describing the repeated interactions of the small system with the heat bath, is given by the sequence $(V_k)_{k\in\N}$ in $\Bb(\Hh_S\otimes\bigotimes_{\N^*}\C^{n+1})$ satisfying
\begin{eqnarray}\label{eq-bath}
\left\{
\begin{array}{lcc}
V_{k+1}=U_{k+1}V_k\\
V_0=I.
\end{array}
\right.
\end{eqnarray}
Note that the operator $U$ can be written as
$$U=\sum_{i,j=0}^n\,U_j^i\otimes a_j^i,$$
where $U^i_j$ are operators on $\Hh_S$. They are the coefficients of
the matrix $(U^i_j)_{0\leq i,j\leq n}$ of $U$ with respect to the basis $\Bb$. Therefore, the equation (\ref{eq-bath}) is written in terms of discrete quantum noises as follows:
\begin{eqnarray*}
\left\{
\begin{array}{lcc}
V_{k+1}=\sum_{i,j=0}^n\,U_j^iV_ka_j^i(k+1)\\
V_0=I.
\end{array}
\right.
\end{eqnarray*}

Next, we give the matrix representation of the operator $U$ with respect to the basis $\Bb$, which will be used. Set $D=(D_{ij})_{1\leq i,j\leq n}$ and consider the matrix $M=(M_{ij})_{1\leq i,j\leq n}$, where $M_{ij}=\delta_{ij}(H_S+\gamma_iI)$. Note that the unitary evolution $U$ can be written as
$$U=e^{-ihH}=\sum_{m\geq0}\frac{(-i)^m}{m!}h^mH^m.$$
Moreover, $H_R$ is a diagonal operator with respect to the basis $\Bb$
$$H_R={\rm diag}\; (\gamma_0,\gamma_1,...,\gamma_n),$$
and the full Hamiltonian $H$ is given by
\begin{eqnarray*}
H=\left(
\begin{array}{cc}
H_S+\gamma_0 I & 0\\
0 & M+\frac{1}{h}D
\end{array}
\right).
\end{eqnarray*}
This implies that
\begin{eqnarray*}
(hH)^2&=&\left(
\begin{array}{cc}
O(h^2) & 0\\
0 & D^2+O(h)
\end{array}
\right).
\end{eqnarray*}
Furthermore, for all $m\geq3$ we get
\begin{eqnarray*}
(hH)^m=\left(
\begin{array}{cc}
o(h^2) & 0\\
0 & D^m+O(h)
\end{array}
\right).
\end{eqnarray*}
Thus, we obtain
\begin{eqnarray}\label{U}
U=\left(
\begin{array}{cc}
I-ih(H_S+\gamma_0I)+O(h^2) & 0\\
0 & I-ihM+(e^{-iD}-I)+O(h)
\end{array}
\right).
\end{eqnarray}
This gives the coefficients $U_j^i$ of the matrix of $U$ with respect to the basis $\Bb$ with precision $O(h)$ and $O(h^2)$.

\subsection{GNS representation}

The aim of this subsection is to describe the GNS representation of the pair $(\C^{n+1},\rho_\beta)$, where $\rho_\beta$ is the thermodynamical state at inverse temperature $\beta$ ( $\beta>0$) of one piece of the heat bath. It is given by
$$\rho_\beta=\frac{1}{Z}e^{-\beta(H_R-\mu N)},$$ 
where
\begin{enumerate}
\item[-] $Z=\Tr(e^{-\beta(H_R-\mu N)}),$
\item[-] $N=\sum_{j=0}^nj|e_j\rangle\langle e_j|$ is the discrete number operator defined on $\C^{n+1},$ 
\item[-] $\mu$ is a scalar, called {\it chemical potential}. 
\end{enumerate}

Note that, with respect to the basis $\Bb$, the density matrix $\rho_\beta$ has the form
 $$\rho_\beta={\rm diag}\;(\beta_0,\beta_1,...,\beta_n),$$
 where
\begin{equation}\label{beta}
\beta_j=\frac{e^{j\mu\beta}e^{-\beta\gamma_j}}{e^{-\beta\gamma_0}+e^{\mu\beta}e^{-\beta\gamma_1}+...+e^{n\mu\beta}e^{-\beta\gamma_n}},
\end{equation}
for all $j\in\{0,1,...,n\}$.

Now, denote by $\widetilde{\Hh}=\Bb(\C^{n+1})$, the algebra of all bounded operators on $\C^{n+1}$ equipped with the scalar product
$$\langle A,B\rangle=\Tr(\rho_\beta A^*B).$$
Hence, the GNS representation of the pair $(\C^{n+1},\rho_\beta)$ is the triple $(\pi,\widetilde{\Hh},\Omega_R)$, where
\begin{enumerate}
\item[-] $\Omega_R=I$,
\item[-] $\pi:\widetilde{\Hh}\longrightarrow \Bb(\widetilde{\Hh})$ such that $\pi(M)A=MA$ for all $M,A\in\widetilde{\Hh}.$
\end{enumerate}

Set
$$\widetilde{U}=\pi(U)$$
and denote by $\widetilde{\Hh}_k$ the $k$-copy of $\widetilde{\Hh}$ in the chain $\bigotimes_{\N^*}\widetilde{\Hh}$. Then, it is easy to check that $\widetilde{U}_k=\pi(U_k)$ acts as $\widetilde{U}$ on $\Hh_S\otimes\widetilde{\Hh}_k$ and the identity elsewhere. Moreover, if we denote by $\widetilde{V}_k=\pi(V_k)$, then it is straightforward to check that the sequence $(\widetilde{V}_k)_{k\in\N}$ in $\Bb(\Hh_S\otimes\bigotimes_{\N^*}\widetilde{\Hh})$ satisfies the following equation: 
\begin{eqnarray}\label{evol}
\left\{
\begin{array}{lcc}
\widetilde{V}_{k+1}=\widetilde{U}_{k+1}\widetilde{V}_k\\
\widetilde{V}_0=I.
\end{array}
\right.
\end{eqnarray}
\section{The atom field}

The space $T\Phi$ given in subsection \ref{subsection} has a
continuous version whose structure we describe below. We refer the interested reader to [At1] for more details.

In what follows, we preserve the same notations as in subsection \ref{subsection} and denote by $\Hh'$ the closed subspace of $\Hh$ generated by vectors $(e_i)_{i\in J}$. The symmetric Fock space constructed over the Hilbert space $L^2(\R_+,\Hh')$ is denoted by $\Phi=\Gamma_S(L^2(\R_+,\Hh'))$ with vacuum vector $\Omega$. The space $\Hh'$ is called the {\it multiplicity space} and $\dim\Hh'$ is called the {\it  multiplicity} of the Fock space $\Gamma_S(L^2(\R_+,\Hh'))$. Now, in order to justify the equality $\Phi=\bigotimes_{\R_+}\Hh$, we introduce the so called {\it Guichardet interpretation} of the Fock space $\Phi$.
\subsection{Guichardet interpretation of the Fock space $\Gamma_S(L^2(\R_+,\Hh))$}
Note that we have the following identification
$$L^2(\R_+,\Hh')= L^2(\R_+\times J,\C),$$
obtained by identifying a vector $f$ in the former space with the function on $\R_+\times J$ defined by $(t,j)\mapsto f_j(t)=\langle v_j/f(t)\rangle$. Therefore, the symmetric Fock space is identified to
$$\bigoplus_{k=0}^\infty L^2_{sym}((\R_+\times J)^k,\C)$$
consisting of vectors $\Psi=(\Psi_k)_{k\geq0}$ such that $\Psi_k\in L^2_{sym}((\R_+\times J)^k,\C)$ and 
$$\|\Psi\|^2_{\Gamma_s(L^2(\R_+,\Hh'))}=\sum_{k\geq0}\frac{1}{k!}\|\Psi_k\|^2_{L^2_{sym}}((\R_+\times J)^k,\C).$$
Denote by $\Sigma_k$ the $k$-standard simplex in $\R^k$. Then, it is straightforward to check that
$$\bigoplus_{k=0}^\infty L^2_{sym}((\R_+\times J)^k,\C)\simeq\bigoplus_{k=0}^\infty L^2(\Sigma_k\times J^k)$$
and 
$$\|\Psi\|^2_{\Gamma_s(L^2(\R_+,\Hh'))}=\sum_{k\geq0}\|\Psi_k\|^2_{L^2(\Sigma_k\times J^k,\C)}.$$
Let $\Pp_k$ be the set of $k$-element $\sigma$ of $\R_+\times J$
$$\sigma=\{(t_1,i_1),...,(t_k,i_k)\},$$
where $t_i\neq t_j$ for all $i\neq j$. It is interesting to note that there is isomorphism from $\Pp_k$ into $\Sigma_k\times J^k$ given by $\{(t_1,i_1),...,(t_k,i_k)\}\longmapsto((t_1,t_2,...,t_k),(i_1,i_2,...,i_k))$ such that $t_1<t_2<...<t_k$. Hence, $\Pp_k$ inherits a measured space structure of  $\Sigma_k\times J^k$.

Set $\Pp_0=\{\emptyset\}$ for which we associate the measure $\delta_\emptyset$ and denote by $d\sigma$ the measure on $\Pp=\cup_k\Pp_k$. Let $\mathcal{F}$ be the associated $\sigma$-field. Then, the Fock space $\Gamma_S(L^2(\R_+,\Hh'))$ is the space $L^2(\Pp,\mathcal{F},d\sigma)$ and the elements of $\Gamma_S(L^2(\R_+,\Hh'))$ are the measurable functions $f'$s from $\Pp$ into $\C$ such that
$$\|f\|^2=\int_{\Pp}|f(\sigma)|^2d\sigma<\infty.$$

In the following, any element $\sigma\in\Pp$ is identified with a family $(\sigma_i)_{1\leq i\leq N}$ of subsets of $\R_+$ such that
$$\sigma_i=\{s\in\R_+,\;(s,i)\in\sigma\}.$$

In order to justify that $\Phi$ is the continuous version of the atom chain $T\Phi$, we need to describe an important representation in the Fock space $\Phi$. For this purpose, we introduce the curve family $\chi_t^i$ defined by
\begin{eqnarray*}
\chi_t^i(\sigma):=
\left\{
\begin{array}{lcc}
\mathds{1}_{[0,t]}(s) &{}&\mbox{ if } \sigma=\{(s,i)\}\\
0 &{}&\mbox{ elsewhere }.
\end{array}
\right.
\end{eqnarray*}
This family satisfies the following:
\begin{enumerate}
\item[-] $\chi_t^i\in\Phi_{(0,t)}=\Gamma_s(L^2((0,t),\Hh')),$ 
\item[-] $\chi_t^i-\chi_s^i\in\Phi_{(s,t)}=\Gamma_s(L^2((s,t),\Hh'))\;\,\mbox{ for all }\,s,t$ such that $\,s\leq t.$
\item[-] $\chi_t^i$ and $\chi_s^j$ are orthogonal elements $\Phi$ for all $i,j$ such that $i\neq j$.
\end{enumerate}
The above properties of the family $\chi^i_t$ allow us to define the {\it Ito integral} in $\Phi$. Consider a family $g=\{g_t^i, \,t\geq0,\,i\in J\}$ of elements in $\Phi$ which satisfies the following:
 \begin{enumerate}
\item[i)] $t\mapsto\|g_t^i\|$ is measurable, for all $i$,
\item[ii)] $g_t^i \in\Phi_{(0,t)}$ for all $t$,
\item[iii)] $\sum_{i\in J}\int_0^\infty \|g_t^i\|^2dt<\infty.$
\end{enumerate}
Such a family is called an {\it Ito integrable} family. Then, if we consider a partition $\{t_j,\;j\in\N\}$ of $\R_+$ with diameter $\delta$ and if we denote by $P_t$ the orthogonal projection on $\Phi_{(0,t)}$, the Ito integral of $g$, $I(g)=\sum_{i\in J}\int^\infty_0g_t^id\chi_t^i$, is the limit in $\Phi$ of 
$$\sum_{i\in J}\sum_{j=0}^\infty\,\frac{1}{t_{j+1}-t_j}\,\int_{t_j}^{t_j+1}\,P_{t_j}g_s^i\,ds\otimes(\chi_{t_j+1}^i-\chi_{t_j}^i),$$
when $\delta$ tends to $0$.
\begin{theo}
The Ito integral $I(g)=\sum_{i\in J}\int^\infty_0g_t^id\chi_t^i$ of an Ito integrable family $g=\{g_t^i, \,t\geq0,\,i\in J\}$ is the element of $\Phi$, given by
\begin{eqnarray*}
I(g)(\sigma)=
\left\{
\begin{array}{lcc}
g_{\vee\sigma}^i(\sigma_-)\;\,&{}& \mbox{ if }\, \vee\sigma\in\sigma_i\\
0  &{}&\mbox{ elsewhere },          
\end{array}
\right.
\end{eqnarray*}
where $\vee\sigma=\sup\big\{t\in\R_+ \mbox{ s.t there exists k which satisfies } (t,k)\in\sigma\big\}$ and\\
$\sigma_-=\sigma\setminus (\vee\sigma,i)\;\mbox{ if }\,(\vee\sigma,i)\in\sigma$.
Moreover, the following isometry formula holds
$$\|I(g)\|^2=\big\|\sum_i\int^\infty_0g_t^id\chi_t^i\big\|^2=\sum_i\int_0^\infty\,\|g_t^i\|^2dt.$$
\end{theo}

Now, consider a family $f=(f^i)_{i\in J}$ of elements of $L^2(\Pp_1)=L^2(\R_+\times J)$. It is obvious that $\{f^i(t)\Omega,\,t\in\R_+,\,i\in J\}$ is an Ito integrable family with Ito integral is given by
$$I(f)=\sum_{i\in J}\int_0^\infty f^i(t)\Omega d\chi_t^i$$ 
and we have
$$I(f)(\sigma)=\left\{
\begin{array}{lcc}
f^i(s)\;\,&{}&\mbox{ if } \,\sigma=\{s\}_i\\
0  &{}&\mbox{ elsewhere }.          
\end{array}
\right.
$$
In the same way, we define the Ito integral of a family  $f\in
L^2(\Pp_k)$ recursively as:
\begin{eqnarray*}
I_k(f)&=&\sum_{i_1,...,i_k=1}^N\int_0^\infty\int_0^{t_k}\,...\,\int_0^{t_2}f_{i_1,...,i_k}(t_1,...,t_k)\Omega\,d\chi_{t_1}^{i_1}\,...\,d\chi_{t_k}^{i_k}\\
    &=&\int_{\Pp_k}f(\sigma)d\chi_{t_t}^{i_1}...d\chi_{t_k}^{i_k}.
\end{eqnarray*}
Moreover, we have
\begin{eqnarray*}
[I_k(f)](\sigma)=
\left\{
\begin{array}{lcc}
f_{i_1,...,i_k}(t_1,...,t_k)\;\, &{}&\mbox{ if }\,\sigma=\{(t_1,i_1)\cup...\cup (t_1,i_k)\}\\
0    &{}&\mbox{ elsewhere }.
\end{array}
\right.
\end{eqnarray*}

Finally, if $f=(f_k)_{k\in\N}\in L^2(\Pp)$, then $I(f)$ is given by
 $$f(\emptyset)\Omega+\sum_{k=1}^\infty I_k(f).$$

The following theorem gives the {\it chaotic representation} of an element $f$ in $\Phi$ (cf [At1]).
\begin{theo}\label{KO}
Every element $f$ of $\Phi$ has a unique chaotic representation 
$$f=\int_{\Pp}\,f(\sigma)d\chi_\sigma$$
which satisfies the isometry formula
$$\|f\|^2=\int_{\Pp}\,|f(\sigma)|^2d\chi_\sigma.$$
\end{theo}
From the above theorem, the space $\Phi$ is interpreted as the continuous version of the space $T\Phi$, where the countable orthonormal basis $\{X_A,\;A\in\Pp_{\N^*,J}\}$ of $T\Phi$ is replaced by the continuous orthonormal basis $\{d\chi_\sigma,\,\sigma\in\Pp\}$ of $\Phi$. 

\subsection{Continuous quantum noises}
The symmetric Fock space
$\Phi=\Gamma_S(L^2(\R_+,\Hh)=\bigotimes_{\R_+}\Hh$ is the natural
space in which we define the annihilation, creation and conservation
operators, which are called {\it continuous quantum noises}. These operators are merely considered as a source of noise, which occurs during the interaction of a quantum system and an exterior system. They are defined by
\begin{eqnarray*}
\lbrack a_i^0(t)f\rbrack(\sigma)&=&\sum_{s\in\sigma_i,\,s\leq t}\,f(\sigma\setminus\{s\}_i),\\
\lbrack a_{0}^i(t)f\rbrack(\sigma)&=&\int_0^t f(\sigma\cup\{s\}_i)\,ds,\\
\lbrack a_j^i(t)f\rbrack(\sigma)&=&\sum_{s\in\sigma_j,\,s\leq t}\,f(\sigma\setminus\{s\}_j\cup\{s\}_i),\\
\lbrack a_i^{0}(t)f\rbrack(\sigma)&=&tf(\sigma).
\end{eqnarray*}
A common domain of these operators is given by
$$\mathcal{D}=\{f\in\Phi,\;\int_\mathcal{P}|\sigma||f(\sigma)|d\sigma<\infty\}.$$
 
The coherent vector $e(f)$ of an element $f\in L^2(\R_+,\Hh')$ is defined by 
$$[e(f)](\sigma)=\prod_{i\in J}\prod_{s\in\sigma_i}f_i(s).$$
In [At2], it is proved that the continuous quantum noises satisfy the following relation
$$\langle e(f),\,a^i_j(t)e(g)\rangle=\int_0^t\bar{f_i}(s)g(s)\,ds\,\langle e(f),\,e(g)\rangle,$$
where $a_0^0(t)=tI$, $h_0(t)=1$,  $t\geq0$ and $h\in
L^2(\R_+,\,\Hh')$. We also have the following table\smallskip\\

\begin{tabular}{|c|c|c|c|}
\hline
          &  $\Omega$ & $d\chi_t^i$ & $d\chi_t^j,\,i\neq j$\\

\hline
$da_i^0(t)$ &  $d\chi_t^i$ &   0     &  0\\
\hline
$da^i_0(t)$ &  $dt I$      &   0  &  0\\
\hline
$da^i_j(t)$ & 0    & $d\chi_t^i$ & 0\\
\hline
\end{tabular} 
\bigskip \\
 
As a corollary of the above table, it is easy to show that the actions of the continuous quantum noises $da_i^j(t),\,i,j\in J\cup\{0\}$ on the element of the orthonormal basis $\{d\chi_\sigma,\,\sigma\in\Pp\}$ of $\Phi$ are similar to the ones of the discrete quantum noises on the elements of the basis $\{X_A,\;A\in\Pp_{\N^*,J}\}$ of $T\Phi$. We then have 
\begin{eqnarray*}
da^i_j(t)d\chi_\sigma&=&d\chi_{\sigma\setminus\{(t,i)\}\cup\{(t,j)\}}\,\mathds{1}_{(t,i)\in\sigma}, \;\mbox{ for all } \,i\neq0,\,j\neq0,\\
da^i_0(t)d\chi_\sigma&=&d\chi_{\sigma\setminus\{(t,i)\}}dt\,\mathds{1}_{(t,i)\in\sigma},\\
da^0_j(t)d\chi_\sigma&=&d\chi_{\sigma\cup\{(t,j)\}}\,\mathds{1}_{(t,0)\in\sigma},\\
da^0_0(t)d\chi_\sigma&=&d\chi_\sigma dt\,\mathds{1}_{\{(n,k)\notin A, \forall k\in J\}}.
\end{eqnarray*}   
\subsection{Quantum Langevin equation}
The quantum Langevin equations or Hudson-Parthasarathy equations play an important role in describing the irreversible evolution of a quantum system in interaction with an exterior system. The ingredients of these equations are quantum noises, which are defined in the previous subsection, and system operators, which control the interaction between the two physical systems.

Let $H,\;L^i_j$ and $S^i_j$ be bounded operators on a separable Hilbert space $\Hh_0$ such that $H=H^*$ and $(S^i_j)_{i,j}$ is unitary. The sum $\sum_kL^{*}_kL_k$ is assumed to be strongly convergent to a bounded operator. Suppose that the operators $L^i_j$ satisfy
\begin{eqnarray}\label{sys}
L^0_0&=&-(iH+\frac{1}{2}\sum_kL^{0*}_kL^0_k),\nonumber\\ 
L_j^0&=&L_j,\nonumber\\
L^i_0&=&-\sum_kL_k^*S_i^k,\\
L_i^j&=&S_i^j-\delta_{ij}I.\nonumber
\end{eqnarray}

Then, we have the following theorem (cf [M],\,[P]).
\begin{theo}
Suppose that the operators $L^i_j$ are given by (\ref{sys}). Then there exists a unique strongly continuous unitary process $U_t$ which satisfies the following stochastic differential equation
$$U_t=I+\sum_{i,j}\int_0^tL_i^jU_sda_i^j(s).$$
\end{theo}
\section{The continuous limit}
In this section we state the main result, which allows us to obtain
the quantum Langevin equations from discrete models. We further study
the convergence to a quantum Langevin equation of the repeated quantum interaction model, describing a small system in interaction with an exterior system, which is introduced in subsection 2.2. Hence, we establish the relation between the time scale $h$, given in (\ref{scale}), and the classical low density limit.

\subsection{Convergence to quantum Langevin equation}
Consider a repeated quantum interaction model of a small system $\Hh_S$ in interaction with an atom chain $\bigotimes_{\N^*}\Hh$, where $\{e_i,\;i\in I\cup\{0\}\}$ is an orthonormal basis of the Hilbert space $\Hh$. Let $U=e^{-ihH}$ be the unitary evolution describing the small system in interaction with a single piece of the chain where $H$ is the associated repeated quantum interaction Hamiltonian. Hence, the discrete evolution equation is given by the sequence $(V_k)_{k\in\N}$ in $\Bb(\Hh_S\otimes\bigotimes_{\N^*}\Hh)$ which satisfies
\begin{eqnarray}\label{equation}
\left\{
\begin{array}{lcc}
V_{k+1}=U_{k+1}V_k\\
V_0=I.
\end{array}
\right.
\end{eqnarray} 

Let $(U_j^i)_{i,j}$ be the matrix of $U$ with respect to the basis $\{e_i,\;i\in I\cup\{0\}\}$. Then, we have the following result (cf [AtP]).
\begin{theo} \label{repeated} Assume that there exist bounded operators $L_i^j,\;i,j\in I\cup\{0\}$ on $\Hh_S$ such that
$$\lim_{h\rightarrow0}\,\frac{U^j_i(h)-\delta_{ij}I}{h^{\varepsilon_{ij}}}=L_i^j,$$
where $\varepsilon^{ij}=\frac{1}{2}(\delta_{i0}+\delta_{0j})$. Assume that the  quantum Langevin equation  
$$
\left\{
\begin{array}{lcc}
dV(t)=\sum_{i,j}L^j_iV(t)da_i^j(t)\\
V(0)=I
\end{array}
\right.
$$
has a unique unitary solution $(V_t)_{t\geq0}$. Then for almost all $t$, the solution $V_{[t/h]}$ of (\ref{equation}) converges strongly to $V_t$, when $h$ tends to $0$.
\end{theo}
\subsection{Low density limit}
All this subsection is devoted to study the continuous limit of
the small system in interaction with a heat bath which is presented in
subsection 2.2. The main assumption, according to which we
suppose that the length of the time interaction between the small system and one piece of the heat bath is related to the chemical potential as follows
$$h^2=e^{\beta\mu}.$$
Therefore, it is clear that $h$ tends to $0$ if and only if the fugacity $e^{\beta\mu}$ tends to $0$, that is the chemical potential $\mu$ converges to $-\infty$.

Now, we give an orthonormal basis of $\widetilde{\Hh}=\Bb(\C^{n+1})$. Put 
$$\nu_k=1-\beta_1-\beta_2-...-\beta_k,\,\mbox{ for all }\, k\in\{1,...,n\}.$$ 
Consider the family $\{X^i_j,\;i,\,j\in\{0,1,...,n\}\}$ such that
\begin{enumerate}
\item[-] $X^0_0=I,$
\item[-] $X^i_j=\frac{1}{\sqrt{\beta_i}}\,a^i_j,\; \mbox{ for all }\, i\neq j,$
\item[-] $X_k^k={\rm diag}\,(\lambda_k^0,\lambda^1_k,...,\lambda_k^{k-1},\lambda_k^k,...,\lambda_k^n),$ 
\end{enumerate}
where
\begin{equation}\label{8}
\lambda_1^0=\lambda_1^2=...=\lambda_1^n=\frac{-\sqrt{\beta_1}}{\sqrt{\nu_1}},\;\,\lambda_1^1=\frac{\sqrt{\nu_1}}{\sqrt{\beta_1}},
\end{equation}
and for all $k\in\{2,...,n\}$, we have
\begin{eqnarray}
\lambda_k^0&=&\lambda_k^{k+1}=...=\lambda_k^{n}= \frac{-\sqrt{\beta_k}}{\sqrt{\nu_{k-1}}\,\sqrt{\nu_k}},\label{9}\\
\lambda_k^k&=&\frac{\sqrt{\nu_k}}{\sqrt{\nu_{k-1}}\,\sqrt{\beta_k}},\label{10}\\
\lambda_k^1&=&\lambda_k^2=...=\lambda_k^{k-1}=0.\label{11}
\end{eqnarray}
Hence, it is straightforward to show that the family $\{X^i_j,\;i,\,j\in\{0,1,...,n\}\}$ is an orthonormal basis  of $\widetilde{\Hh}$,
equipped with the scalar product $\langle A,B\rangle=\Tr(\rho_\beta A^*B)$. 

Note that in [AtJ], the authors only need to give explicitly the
elements $X^0_0$ and $X^i_j$, $i\neq j$, because the vectors $X^k_k$
do not contribute in the proof of their main theorem. However, this is not the
case here, so we have computed explicitly the elements $X_k^k$ which
play an important role in the proof of our result.

As a consequence of the relations (\ref{beta}) and (\ref{8})--(\ref{11}), we prove the following.
\begin{lem}\label{lem}
The following hold:
\begin{eqnarray*}
\beta_0\lambda_1^0&=&O(h),\\
\beta_0\lambda_i^0&=&o(h),\;\mbox{ for all }\,i\geq2,\\
\beta_1\lambda_1^1&=&O(h),\\
\beta_k\lambda_i^k&=&o(h),\;\mbox{ for all }\,i\geq1,\; k\geq1\;\mbox{ such that }\,(i,k)\neq(1,1),\\
\beta_k\lambda_i^k\lambda_j^k&=&o(h), \;\mbox{ for all }\,i,j\geq1\, \mbox{ such that }\,i\neq j,\\
\beta_k(\lambda_i^{k})^2&=&o(h),\;\mbox{ for all }\,i,k\geq1\, \mbox{ such that }\,i\neq k,\\
\lim_{h\rightarrow0}\beta_k(\lambda_k^{k})^2&=&1,\;\mbox{ for all }\,k\geq1.
\end{eqnarray*}
\end{lem}

In order to study the continuous limit of the discrete solution of equation (\ref{evol}), we consider the matrix representation $(\widetilde{U}_{k,l}^{i,j})_{i,j,k,l\in\{0,1,...,n\}}$ of $\widetilde{U}$ with respect to the basis $\{X^i_j,\;i,\,j=0,1,...,n\}$. Therefore, it is obvious that
\begin{equation}\label{wide}
\widetilde{U}_{k,l}^{i,j}=\Tr_{\widetilde{\Hh}}(\rho_\beta(X^k_l)^*UX^i_j).
\end{equation}

In the sequel, we suppose that the matrix of the unitary operator $e^{-iD}$, defined by (3), with respect to the basis $\mathcal{B}$, is given by
$$e^{-iD}=(S_l^k)_{1\leq l,k\leq n}.$$

Now, we prove the following.
\begin{theo}
The solution $\widetilde{V}_{[t/h]}$ of (\ref{evol}) converges strongly, when $h$ tends to $0$, to the unitary solution of the quantum Langevin equation
\begin{eqnarray}\label{eq}
d\widetilde{V}_t=&-&i(H_S+\gamma_0I)\widetilde{V}_tdt\nonumber\\
     &+&\sum_{j,k=1}^N\,(S^j_k-\delta_{jk}I)\widetilde{V}_t\Bigl(\sum_{i=1}^nda_{i,k}^{i,j}(t)\Bigr),
\end{eqnarray}
with initial condition $\widetilde{V}_0=I$ and where $da_{l,k}^{i,j}(t),\;i,j,k,l=0,1,...,n$ are the associated quantum noises of the symmetric Fock space $\Gamma_s(L^2(\R_+,\,\C^{(n+1)^2-1}))$ with respect to the basis $\{X^i_j,\;i,\,j=0,1,...,n\}$.
\end{theo}\label{true}
\begin{proof}From relation (\ref{wide}) we have 
$$\widetilde{U}_{0,0}^{0,0}=\Tr_{\widetilde{\Hh}}(\rho_\beta U).$$
Hence, by using (\ref{U}) we get
\begin{eqnarray*}
\widetilde{U}_{0,0}^{0,0}=&{}&\!\!\!\beta_0(I-ihH_S+\gamma_0hI)+\beta_1(I-ihH_S)\\
                         &+&\beta_2(I-ihH_S)+...+\beta_n(I-ihH_S)+O(h^2).
\end{eqnarray*}
This gives
\begin{equation}\label{U0000}
\widetilde{U}_{0,0}^{0,0}=I-ih(H_S+\beta_0\gamma_0I)+o(h).
\end{equation}

Now, for all $i,j\in\{0,1,...,n\}$ such that $(i,j)\neq(0,0)$ and $i\neq j$, we have
\begin{eqnarray*}
\widetilde{U}_{0,0}^{i,j}&=&\Tr_{\widetilde{\Hh}}(\rho_\beta UX^i_j)\\
                       &=&\frac{1}{\sqrt{\beta_i}}\langle e_i,\,\rho_\beta Ue_j\rangle\\
                       &=& \sqrt{\beta_i}\langle e_i,\,Ue_j\rangle.
\end{eqnarray*}
Therefore, we distinguish the following two cases:
\begin{enumerate}
\item[-] If $i=0$ or $j=0$, then we have
\begin{equation}\label{15}
\widetilde{U}_{0,0}^{i,0}=\widetilde{U}_{0,0}^{0,j}=0.
\end{equation}
\item[-] If $i\neq0$ and $j\neq0$, then we have
\begin{equation}\label{16}
\widetilde{U}_{0,0}^{i,j}=
\left\{
\begin{array}{lcc}
O(h) &  \;\mbox{ if }\;i=1\\
o(h) &  \;\mbox{ if }\;i\neq1 .
\end{array}
\right.
\end{equation}
\end{enumerate}

In the same way, we prove that
\begin{equation}\label{17}
\widetilde{U}_{0,l}^{0,0}=\widetilde{U}_{k,0}^{0,0}=0,\;\forall k,l\in\{1,...,n\},
\end{equation}
and for all $k,l\in\{1,...,n\}\,\mbox{ such that }\,k\neq l$, we get 
\begin{equation}\label{18}
\widetilde{U}_{k,l}^{0,0}=
\left\{
\begin{array}{lcc}
O(h) & \;\mbox{ if }\,k=1\\
o(h) & \;\mbox{ if }\,k\neq1.
\end{array}
\right.
\end{equation}

It is worthwhile to note that for all $i\neq j,\;k\neq l$ and $i,j,k,l\in\{1,...,n\}$, we have
$$\widetilde{U}_{k,l}^{i,j}=\delta_{ik}\langle e_l,Ue_j\rangle,$$
from which follows that
\begin{equation}\label{19}
\widetilde{U}_{k,l}^{i,j}=0,\;\forall i\neq k
\end{equation}
and
\begin{equation}\label{20}
\widetilde{U}_{i,l}^{i,j}=\langle e_l,Ue_j\rangle=S_l^j+O(h).
\end{equation}

For all $i\in\{1,...,n\}$, the coefficient $\widetilde{U}_{0,0}^{i,i}$ is given by
\begin{eqnarray*}
\widetilde{U}_{0,0}^{i,i}&=&\Tr_{\widetilde{\Hh}}(\rho_\beta UX^i_i)\\
                        &=&\sum_j\beta_j\lambda_i^j\langle e_j,\,Ue_j\rangle.
\end{eqnarray*}
Note that from Lemma \ref{lem}, we have the following
\begin{eqnarray*}
\beta_0\lambda_1^0&=&O(h),\;\beta_1\lambda_1^1=O(h),\\
\beta_0\lambda_i^0&=&o(h),\;\mbox{ for all }i\geq2,\\
\beta_j\lambda_i^j&=&o(h),\;\mbox{ for all }i\geq1,\,j\geq1\;\mbox{ such that }(i,j)\neq(1,1).
\end{eqnarray*}
Hence, we obtain
\begin{equation}\label{21}
 \widetilde{U}_{0,0}^{i,i}=
\left\{
\begin{array}{lcc}
O(h) & \mbox { if } \,i=1\\
o(h) & \mbox { if } \,i\geq2.
\end{array}
\right.
\end{equation}
Similar as a above, we show that
\begin{equation}\label{22}
\widetilde{U}_{k,k}^{0,0}=
\left\{
\begin{array}{lcc}
O(h) & \mbox{ if }\,k=1\\
o(h) & \mbox{ if }\, k\geq2.
\end{array}
\right.
\end{equation}

Now, for all $i,\;k\;l\in\{1,...,n\}$ such that $k\neq l$, we have
\begin{eqnarray*}
\widetilde{U}_{k,l}^{i,i}&=&\Tr_{\widetilde{\Hh}}(\rho_\beta(X^k_l)^*UX^i_i)\\
                        &=&\sum_j\,\frac{\beta_j}{\sqrt{\beta_k}}\lambda_i^j\delta_{kj}\langle e_l,Ue_j\rangle\\
                        &=&\sqrt{\beta_k}\lambda_i^k\langle e_l,Ue_k\rangle.
\end{eqnarray*}
Note that $\sqrt{\beta_k}=o(h),\;\forall k\in\{1,...,n\}$. Therefore, we get
\begin{equation}\label{23}
\widetilde{U}_{k,l}^{i,i}=o(h),\;\forall\,i\neq k
\end{equation}
and
\begin{equation}\label{24}
\widetilde{U}_{i,l}^{i,i}=\frac{\sqrt{\nu_k}}{\sqrt{\nu_{k-1}}}\,S^i_l+O(h).
\end{equation}
In the same way, we prove that for all $i,\,j,\,k\in\{1,...,n\}$ such that $i\neq j$, we have
\begin{equation}\label{25}
\widetilde{U}^{i,j}_{k,k}=o(h),\;\forall i\neq k
\end{equation}
and
\begin{equation}\label{26}
\widetilde{U}^{i,j}_{i,i}=\frac{\sqrt{\nu_k}}{\sqrt{\nu_{k-1}}}\,S^j_i+O(h).
\end{equation}
Let $i,k\in\{1,...,n\}$. We then have
\begin{eqnarray*}
\widetilde{U}_{k,k}^{i,i}&=&\Tr_{\widetilde{\Hh}}(\rho_\beta X^k_kUX^i_i)\\
                        &=&\sum_j\,\beta_j\lambda_k^j\lambda_i^j\langle e_j,Ue_j\rangle\\
                        &=&\sum_{j=1}^N\beta_j\lambda_k^j\lambda_i^j(S^j_j-I)+\sum_{j=1}^n\beta_j\lambda_k^j\lambda_i^j\,I+\beta_0\lambda_k^0\lambda_i^0I+o(h).
\end{eqnarray*}
Note that
$$\langle X^i_i,X^k_k\rangle=\sum_{j=0}^n\beta_j\lambda_k^j\lambda_i^j=\delta_{ik}.$$
This implies that
$$\widetilde{U}_{k,k}^{i,i}=\delta_{ik}I+\sum_{j=1}^N\beta_j\lambda_k^j\lambda_i^j(S^j_j-I)+o(h).$$
Thus, we obtain
\begin{equation}\label{27}
\widetilde{U}_{k,k}^{i,i}=o(h),\;\forall\,i\neq k
\end{equation}
and 
\begin{equation}\label{28}
\widetilde{U}_{i,i}^{i,i}=I+\beta_i(\lambda_i^{i})^2(S^i_i-I)+o(h).
\end{equation}

Now, in order to apply Theorem \ref{repeated}, we will compute the following limits
$$s-\lim_{h\rightarrow0}\frac{\widetilde{U}_{k,l}^{i,j}-\delta_{(i,j),(k,l)}I}{h^{\varepsilon_{k,l}^{i,j}}},$$
where $\varepsilon_{0,0}^{0,0}=1,\;\varepsilon_{k,l}^{0,0}=\varepsilon_{0,0}^{i,j}=1/2$ and $\varepsilon_{k,l}^{i,j}=0$. 

Note that from equalities (\ref{15}), (\ref{17}) and (\ref{19}), we have $\widetilde{U}_{0,0}^{i,0}=\widetilde{U}_{0,0}^{0,j}=0$ for all $i,j\in\{1,...,n\}$ and $\widetilde{U}_{k,l}^{i,j}=0$ for all $i,j,k,l\in\{1,...,n\}$ such that $i\neq k$. Moreover, the equalities (\ref{16}), (\ref{18}), (\ref{21}) and (\ref{22}) imply that
$$\lim_{h\rightarrow0}\frac{\widetilde{U}_{k,l}^{0,0}}{\sqrt{h}}=\lim_{h\rightarrow0}\frac{\widetilde{U}_{0,0}^{i,j}}{\sqrt{h}}=0,$$
 $\forall i,j,k,l\in\{1,...,n\}$.

By using relations (\ref{20}), (\ref{24}) and (\ref{26},) we have
$$\lim_{h\rightarrow0}\widetilde{U}_{i,l}^{i,j}=S_l^j,$$
for all $ i,j,l\in\{1,...,n\}$ such that $j\neq l$. Furthermore, by taking into account equalities (\ref{23}), (\ref{25}) and (\ref{27}), we get
$$\lim_{h\rightarrow0}\widetilde{U}_{k,l}^{i,i}=\lim_{h\rightarrow0}\widetilde{U}_{k,l}^{i,j}=\lim_{h\rightarrow0}\widetilde{U}_{k,k}^{i,i}=0,$$
 $\forall i,j,k,l\in\{1,...,n\}$ such that $i\neq k$.

The equality (\ref{U0000}) implies that
$$\lim_{h\rightarrow0}\frac{\widetilde{U}_{0,0}^{0,0}-I}{h}=-i(H_S+\gamma_0I).$$
Finally, from (\ref{20}) and (\ref{28}), we have
$$\lim_{h\rightarrow0}(\widetilde{U}_{i,j}^{i,j}-I)=S^j_j-I,\;\forall i,j\in\{1,...,n\}.$$
Hence, by using Theorem \ref{repeated}, the solution of equation (\ref{evol}) converges strongly to the unitary solution of the quantum stochastic differential equation
\begin{eqnarray*}
d\widetilde{V}_t&=&-i(H_S+\gamma_0I)\,\widetilde{V}_tdt\\
     &+&\sum_{\scriptsize \begin{array}{ll} j,k=1\\ j\neq k\end{array}}^nS_k^j\widetilde{V}_t\sum_{i=1}^nda_{i,k}^{i,j}(t)\\
     &+&\sum_{j=1}^n\,(S^j_j-I)\,\widetilde{V}_t\sum_{i=1}^nda_{i,j}^{i,j}(t)\\
    &=& -i(H_S+\gamma_0I)\,\widetilde{V}_tdt\\
    &+&\sum_{j,k}^n(S_k^j-\delta_{j,k}I)\widetilde{V}_t\sum_{i=1}^nda_{i,k}^{i,j}(t),
\end{eqnarray*}
with initial condition $\widetilde{V}_0=I$. This completes the proof.
\end{proof}

Setting 
\begin{equation}\label{H-P}
dA_k^j(t)=\sum_{i=1}^nda_{i,k}^{i,j}(t),\;\mbox{ for all } j,k\geq1,
\end{equation}
 By the following proposition, we prove that the noises defined by (\ref{H-P}) satisfy the usual Hudson-Parthasarathy {\it Ito table} (cf \cite{HP}).

\begin{prop}
For all $i,\,j\geq1$ we have
$$dA_k^j(t)dA_l^m(t)=\delta_{jl}dA_k^m(t).$$
\end{prop}
\begin{proof}
We have
\begin{eqnarray*}
dA_k^j(t)dA_l^m(t)&=&\Big(\sum_{i=1}^nda_{i,k}^{i,j}(t)\Big)\Big(\sum_{i=1}^nda_{i,l}^{i,m}(t)\Big)\\
 &=&\sum_{i_1,i_2=1}^n\,da_{i_1,k}^{i_1,j}(t)da_{i_2,l}^{i_2,m}(t)\\
 &=&\sum_{i_1,i_2=1}^n\,\delta_{(i_1,j),(i_2,l)}da_{i_1,k}^{i_2,m}(t)\\
 &=&\delta_{jl}\sum_{i=1}^nda_{i,k}^{i,m}(t)\\
 &=&\delta_{jl}dA_k^m(t).
\end{eqnarray*}
\end{proof}

Therefore, from the above proposition the quantum noises $dA_k^j(t)$ satisfy the well-known Ito formula. Moreover, the equation (\ref{eq}) can be identified to a quantum Langevin equation defined on a Fock space with multiplicity $n^2$. Hence, on the Fock space $\Hh_S\otimes\Gamma_s(L^2(\R_+,\C^{n^2}))$, this equation is written in terms of quantum noises $dA_k^j(t),\,i,j\geq1$ as 
\begin{eqnarray*}
d\widetilde{V}_t&=&-i(H_S+\gamma_0I)\,\widetilde{V}_tdt\\
                &+&\sum_{j,k}^n(S_k^j-\delta_{jk}I)\widetilde{V}_tdA_k^j(t).     
\end{eqnarray*}

It is worth noticing that in [AtJ], the repeated quantum interaction
Hamiltonian studied by the authors is composed of a free part and a
dipole Hamiltonian interaction with time scale $\sqrt{h}$. They prove
that, in the continuous limit, we get a quantum diffusion equation,
where new noises appear, called {\it the thermal noises}. Moreover,
these noises satisfy a specified commutation relation, which depends
on the temperature $\beta$. In our paper, after
taking the continuous limit of the repeated quantum interaction model
with time scale $h$, we obtain a quantum Langevin equation directed by Poisson processes. Furthermore, the noises, which appear in this equation, have the same properties as the canonical continuous quantum noises.

{\small


\begin{thebibliography}{99}
\bibitem[At1]{} S. Attal: {Extensions of the quantum stochastic calculus, {\it Quantum Probability Communications XI, 1-38, World Scientific (2003).}} 
\bibitem[At2]{} S. Attal: {\it Quantum noises theory, {submitted}}.
\bibitem[AFL]{} L. Accardi, A. Frigerio, Y.G. Lu: {Weak coupling limit as a quantum functional central limit theorem}. {\it Com. Math. Phys. 131, 537-570 (1990).}
\bibitem[APV]{} L. Accardi, A.N. Pechen, I. V. Volovich: {A Stochastic Golden Rule and Quantum Langevin Equation For the Low Density Limit. {\it Infinite Dimensional Analysis, Quantum Probability and Related Topics, Vol. 6, No. 3 (2003) 431-453.}}
\bibitem[AtP]{}  S. Attal, Y. Pautrat: {From Repeated to Continuous Quantum Interactions. {\it Annales Henri Poincar\'e (Physique Th\'eorique) 7 (2006), p. 59-104.}}
\bibitem[AtJ]{}  S. Attal, A. Joye: {The Langevin equation for a quantum heat bath. {\it Journal of Functional Analysis, to appear}}.
\bibitem[Dav]{} E.B. Davies: {Markovian Master equations. {\it Comm. Math. Phys. 39 (1974), 91-110.}}
\bibitem[D1]{} A. Dhahri: {Markovian Properties of the spin-boson model. {\it S\'eminaire de Probabilit\'es, to appaear}}.
\bibitem[D2]{} A. Dhahri: {A Lindblad model for a spin chain coupled to thermal heat baths. {\it preprint}}.
\bibitem[G]{} A. Guichardet: {Symmetric Hilbert spaces and related Topics. {\it Lectures Notes in Mathematics 261, Springer Verlag (1972)}}.
\bibitem[HP]{HP} R.L. Hudson, K.R. Parthasaraty: {Quantum It\^o's formula and stochastic evolutions. {\it Communications in Mathematical Physics 93 (1984), p. 301-323}}.
\bibitem[M]{} P. A. Meyer: {{\it Quantum Probability for Probabilists }, Second edition. Lect Not. Math. 1538, Berlin: Springer-Verlag 1995.}
\bibitem[P]{} K. R. Parthasarathy: {{\it An Introduction to Quantum Stochastic Calculus}, Birkh\"auser Verlag: Basel. Boston. Berlin.}
\bibitem[Pe]{} A. N. Pechen: {Quantum stochastic equation for test particle interacting with dilute bose gas, {\it  J. Math. Phys. 45 (2004), no. 1, 400-417}}.
\end{thebibliography}
\end{document}